\documentclass[11pt,a4paper]{article}

\usepackage{amsmath}
\usepackage{amsthm}
\usepackage{amssymb}
\usepackage{amscd}
\usepackage[all]{xy}

\title{Effective Non-vanishing for Algebraic Surfaces \\
in Positive Characteristic\footnote{{\it 2000 Mathematics Subject 
Classification} Primary 14J26; Secondary 14E30.}}
\author{Qihong Xie}
\date{\small Dedicated to Professor Musheng Yao on the occasion of 
his sixtieth birthday}
\theoremstyle{plain}
\newtheorem{prop}{Proposition}[section]
\newtheorem{lem}[prop]{Lemma}
\newtheorem{thm}[prop]{Theorem}

\newtheorem{conj}[prop]{Conjecture}
\newtheorem{prob}[prop]{Problem}

\newtheorem*{case}{Case}

\theoremstyle{definition}
\newtheorem{defn}[prop]{Definition}
\newtheorem*{ack}{Acknowledgements}

\theoremstyle{remark}
\newtheorem{rem}[prop]{Remark}
\newtheorem{ex}[prop]{Example}

\newcommand{\Q}{\mathbb Q}
\newcommand{\R}{\mathbb R}

\newcommand{\F}{\mathbb F}
\newcommand{\PP}{\mathbb P}
\newcommand{\OO}{\mathcal O}

\newcommand{\EE}{\mathcal E}

\newcommand{\LL}{\mathcal L}

\newcommand{\FF}{\mathcal F}
\newcommand{\BB}{\mathcal B}
\newcommand{\cP}{\mathcal P}
\newcommand{\Pic}{\mathop{\rm Pic}\nolimits}
\newcommand{\Alb}{\mathop{\rm Alb}\nolimits}
\newcommand{\Supp}{\mathop{\rm Supp}\nolimits}

\newcommand{\ch}{\mathop{\rm char}\nolimits}
\newcommand{\Exc}{\mathop{\rm Exc}\nolimits}

\newcommand{\kolv}{\hbox{\rm Koll\'ar vanishing}}
\newcommand{\spa}{\hbox{\rm semipositivity}}
\newcommand{\lsp}{\hbox{\rm log. semipositivity}}
\newcommand{\en}{\hbox{\rm effective non-vanishing}}
\newcommand{\cc}{\hbox{\rm cyclic cover}}
\newcommand{\kodv}{\hbox{\rm Kodaira vanishing}}
\newcommand{\kvv}{\hbox{\rm Kawamata-Viehweg vanishing}}

\begin{document}
\maketitle

\begin{abstract}
We give a theorem on the effective non-vanishing problem 
for algebraic surfaces in positive characteristic. 
For the Kawamata-Viehweg vanishing, the logarithmic Koll\'ar 
vanishing and the logarithmic semipositivity, we give their 
counterexamples on ruled surfaces in positive characteristic.
\end{abstract}

\setcounter{section}{0}
\section{Introduction}\label{D:S1}

In this paper, we shall consider the following effective 
non-vanishing problem.

\begin{prob}\label{D:1.1}
Let $X$ be a normal proper algebraic variety over an algebraically 
closed field $k$, and $B=\sum b_iB_i$ an effective $\R$-divisor on $X$ 
such that $(X,B)$ is Kawamata log terminal. Let $D$ be a nef Cartier 
divisor on $X$ such that $H=D-(K_X+B)$ is nef and big. 
Find the smallest positive integer $m$ such that $H^0(X,mD)\neq 0$.
\end{prob}

In this problem, we may require the smallest positive integer $m$ 
is universal in the sense that it depends only on the dimension 
of $X$. Furthermore, Ambro and Kawamata have conjectured that, 
if the characteristic of $k$ is zero then $m$ is equal to one, 
which is called the effective non-vanishing conjecture 
(cf.\ \cite{am, ka00}).

\begin{conj}[Effective Non-vanishing]\label{D:1.2}
With the same assumptions as in Problem \ref{D:1.1}, assume further 
that $\ch(k)=0$. Then $H^0(X,D)\neq 0$ holds.
\end{conj}

For the convenience of the reader, we give some necessary definitions. 

\begin{defn}\label{D:1.3}
Let $X$ be a normal proper algebraic variety over an algebraically 
closed field $k$, and $B=\sum b_iB_i$ an effective $\R$-divisor on $X$. 
The pair $(X,B)$ is said to be Kawamata log terminal (KLT, for short), 
or to have Kawamata log terminal singularities, 
if the following conditions hold:

(1) $K_X+B$ is $\R$-Cartier, i.e.\ $K_X+B$ is an $\R$-linear combination 
of Cartier divisors;

(2) For any birational morphism $f:Y\rightarrow X$, we may write 
$K_Y+B_Y\equiv f^*(K_X+B)$, where $\equiv$ means numerical equivalence, 
and $B_Y=\sum a_iE_i$ is an $\R$-divisor on $Y$. 
Then $a_i<1$ hold for all $i$.

Firstly, it follows from (2) that $[B]=0$, i.e.\ $b_i<1$ for all $i$, 
where $[B]=\sum [b_i] B_i$ is the round-down of $B$.

Secondly, this definition is characteristic free, hence it makes sense in 
positive characteristic as well as Problem \ref{D:1.1}.

Thirdly, provided that $\ch(k)=0$ or $\dim X\leq 2$, then $X$ admits a 
log resolution, i.e.\ there exists a desingularization $f:Y\rightarrow X$ 
from a nonsingular variety $Y$, such that the union of the strict 
transform $f_*^{-1}B$ of $B$ and the exceptional locus $\Exc(f)$ of $f$ 
has simple normal crossing support. In this time, condition (2) holds 
for all birational morphisms is equivalent to that it holds for a log 
resolution of $X$. Note that the existence of resolution of singularities 
in positive characteristic is conjectural for higher dimensions.

Let us mention a simple example of KLT pair. Let $X$ be a nonsingular 
variety and $B$ an effective $\Q$-divisor on $X$ such that 
$[B]=0$ and $\Supp(B)$ is simple normal crossing. 
Then $(X,B)$ is KLT.

Assume that $B=0$ and that $K_X$ is $\Q$-Cartier. 
We write $K_Y\equiv f^*K_X+\sum a_iE_i$, where $E_i$ are all 
exceptional divisors of $f$. Then $X$ is said to be terminal, 
or have terminal singularities if $a_i>0$ for all $i$ and for all 
birational morphisms $f:Y\rightarrow X$. For instance, when $\dim X=2$, 
$X$ is terminal is equivalent to that $X$ is nonsingular. 
Reid and Mori gave the classification of all 3-dimensional terminal 
singularities.

Similarly, we can give the definitions of other types of singularities, 
such as canonical, purely log terminal, divisorial log terminal 
and log canonical. We refer the reader to \cite{km} for more details.

Let $D$ be an $\R$-Cartier divisor on $X$. $D$ is said to be nef, 
if $D.C\geq 0$ holds for any irreducible proper curve $C$ in $X$. 
$D$ is said to be nef and big, if $D$ is nef and if the top 
self-intersection $D^n>0$ holds for $n=\dim X$.
\end{defn}

The assumptions in Problem \ref{D:1.1} are standard from the viewpoint 
of the minimal model theory. We have the following celebrated theorems 
(cf.\ \cite[Theorems 1-2-5,2-1-1,3-1-1]{kmm}):

\begin{thm}[Kawamata-Viehweg Vanishing]\label{D:1.4}
With the same assumptions as in Conjecture \ref{D:1.2}. 
Then $H^i(X,D)=0$ holds for any $i>0$.
\end{thm}

\begin{thm}[Non-vanishing and Base Point Free]\label{D:1.5}
With the same assumptions as in Conjecture \ref{D:1.2}. 
Then for any $m\gg 0$, $H^0(X,mD)\neq 0$ holds and the linear system 
$|mD|$ is base point free.
\end{thm}

Clearly, Problem \ref{D:1.1} is just to find the minimum of 
$m$ in the non-vanishing theorem. From the Kawamata-Viehweg 
vanishing theorem, it follows that in Conjecture \ref{D:1.2}, 
$H^0(X,D)\neq 0$ is equivalent to $\chi(X,D)\neq 0$, which shows 
this conjecture is indeed a topological problem in some sense.

We digress to give the history behind Conjecture \ref{D:1.2}, 
and to convince the reader that this conjecture 
is closely related to the minimal model theory.

A nonsingular variety $X$ of dimension $n$ is called a Fano $n$-fold, 
if the anticanonical divisor $-K_X$ is ample. The Fano index of $X$ is, 
by definition, the greatest positive integer $r$ such that 
$-K_X=rD$ for some integral divisor $D$ on $X$. 
As is well-known, the classification of Fano $n$-folds is one of 
the most important problems in algebraic geometry, not only because 
it is interesting in its own right, but also because Fano variety is 
a kind of outcomes when we run the minimal model program 
for nonsingular varieties.

We may consider the following problem, whose first part is 
a very special case of Conjecture \ref{D:1.2}.

\begin{prob}\label{D:1.6}
Let $X$ be a Fano $n$-fold, $r$ the Fano index of $X$, and $D$ an 
integral divisor such that $-K_X=rD$. 
Do the following problems have affirmative answers?

(1) $H^0(X,D)\neq 0$;

(2) The general member of $|D|$ is a nonsingular $(n-1)$-fold.
\end{prob}

Iskovskikh and Shokurov first studied Fano 3-folds in characteristic 
zero. Iskovskikh classified Fano 3-folds of the first kind (namely, 
Fano 3-folds $X$ with the second Betti number $b_2(X)=1$), under the 
assumption that Problem \ref{D:1.6} is true. Shokurov proved that 
Problem \ref{D:1.6} is indeed true for all Fano 3-folds of the first 
kind, and thereby validated Iskovskikh's classification result of 
Fano 3-folds of the first kind. Mori and Mukai classified all 
Fano 3-folds $X$ with $b_2(X)\geq 2$ by virtue of the extremal 
ray theory. Later, Fujita (case $r=n-1$) and Mukai (case $r=n-2$) 
offered an idea to generalize the Iskovskikh and Shokurov's framework 
of Fano 3-folds to that of Fano $n$-folds with $n\geq 4$, and Mukai 
obtained the classification of Fano 4-folds with $r=2$ provided that 
Problem \ref{D:1.6} is true. Thus it is so clear that Problem 
\ref{D:1.6} is a basis to the classification of Fano $n$-folds.

When running the minimal model program started from a nonsingular 
variety $X$ of dimension $n\geq 3$, we have to consider the singularities. 
It turns out that the category of terminal varieties is suitable for 
running the minimal model program. Namely, for any nonsingular variety 
$X$, by virtue of extremal divisorial contractions or flips, 
finally we can obtain a minimal model or a Mori fiber space which 
belongs to this category and is birational to $X$. 
This statement is called the minimal model conjecture, 
and was already proved for dimension $n=3,4$ and $\ch(k)=0$. 
When considering the pair $(X,B)$ with $B$ a suitable $\Q$-divisor on 
a nonsingular variety $X$, the corresponding statement is called 
the log minimal model conjecture, and a suitable category is the 
category of varieties with KLT singularities.

From the viewpoint of the (log) minimal model program, we should consider 
a Fano variety with suitable singularities, since it appears 
as the general fiber of some (log) Mori fiber space.

\begin{defn}\label{D:1.7}
Let $X$ be a normal proper variety of dimension $n$. $X$ is called a 
terminal $\Q$-Fano $n$-fold, if $X$ is terminal and $-K_X$ is ample. 

Let $B$ be an effective $\Q$-divisor on $X$. $(X,B)$ is called a 
KLT $\Q$-Fano pair, if $(X,B)$ is KLT and $-(K_X+B)$ is ample.

The Fano index of a terminal $\Q$-Fano $n$-fold $X$ is the greatest 
rational number $r$ such that $-K_X\sim_\Q rD$ for some Cartier divisor 
$D$ on $X$. The Fano index of a KLT $\Q$-Fano pair $(X,B)$ is the 
greatest rational number $r$ such that $-(K_X+B)\sim_\Q rD$ for some 
Cartier divisor $D$ on $X$.
\end{defn}

The classification of terminal $\Q$-Fano $n$-folds is more difficult 
than that of nonsingular Fano $n$-folds for $n\geq 3$. It seems impossible 
to classify all KLT $\Q$-Fano pairs. However, we can consider the similar 
one to Problem \ref{D:1.6}, whose first part is also a special case of 
Conjecture \ref{D:1.2} (we omit the terminal $\Q$-Fano case below).

\begin{prob}\label{D:1.8}
Let $(X,B)$ be a KLT $\Q$-Fano pair of dimension $n$, 
$r$ the Fano index of $(X,B)$, 
and $D$ a Cartier divisor such that $-(K_X+B)\sim_\Q rD$. 
Do the following problems have affirmative answers?

(1) $H^0(X,D)\neq 0$;

(2) Let $X'\in |D|$ be the general member. Then $(X',B|_{X'})$ has 
KLT singularities.
\end{prob}

The following theorem due to Ambro, gave a partial answer to 
Problem \ref{D:1.8} (cf.\ \cite[Main Theorem]{am}).

\begin{thm}\label{D:1.9}
With the same assumptions as in Problem \ref{D:1.8}, assume further 
that $r>n-3$ and $\ch(k)=0$. Then Problem \ref{D:1.8} is true.
\end{thm}

Note that if $n\leq 3$ then the assumption $r>n-3$ is trivial. Thus 
Problem \ref{D:1.8} is true for all KLT $\Q$-Fano pair of dimension 
$n\leq 3$. When $n=4$, $B=0$ and $X$ has only Gorenstein canonical 
singularities, Kawamata dealt with the case $r=1$ by showing that 
$H^0(X,D)\neq 0$ holds and that the general member $X'\in |D|$ has also 
only Gorenstein canonical singularities (cf.\ \cite[Theorem 5.2]{ka00}).

Let us return to the argument of Conjecture \ref{D:1.2}. 
It is easily verified in the curve case by using the 
Riemann-Roch theorem. The surface case, which is absolutely 
nontrivial, was proved by Kawamata (cf.\ \cite[Theorem 3.1]{ka00}), 
by means of the following so-called logarithmic semipositivity 
theorem (we omit its general statement and only give a special case 
where the base space is 1-dimensional, cf.\ \cite[Theorem 1.2]{ka00}). 

\begin{thm}[Logarithmic Semipositivity]\label{D:1.10}
Let $X$ be a normal proper variety over an algebraically closed field 
$k$ with $\ch(k)=0$, and $B$ an effective $\Q$-divisor on $X$ such that 
$(X,B)$ is KLT. Let $f:X\rightarrow C$ be a surjective morphism to a 
smooth curve $C$. Let $D$ be a Cartier divisor on $X$ such that 
$D\sim_\Q K_{X/C}+B$. Then $f_*\OO_X(D)$ is a semipositive locally 
free sheaf on $C$.

A locally free sheaf $\EE$ on $C$ is said to semipositive, 
if for any morphism $g:C'\rightarrow C$ from a smooth curve 
$C'$ to $C$, and for any quotient line bundle $\LL$ of $g^*\EE$ 
on $C'$, we have $\deg\LL\geq 0$ holds.
\end{thm}

For the higher dimensional cases, the effective non-vanishing 
conjecture is still open, and only a few results are known. We list 
them in the following remark.

\begin{rem}\label{D:1.11}
With the same assumptions as in Conjecture \ref{D:1.2}.

(1) If $(X,B)$ is assumed to be log canonical, but not KLT, 
then we can reduce this case to the KLT case of lower dimension, 
by means of Kawamata's subadjunction theorem and the Nadel vanishing 
theorem (cf.\ \cite[Appendix]{am}). So we only need to treat 
the KLT case from the beginning.

(2) If the irregularity $q(X):=h^1(X,\OO_X)>0$, then we can 
reduce this case to the lower dimensional case, by virtue of 
the Fourier-Mukai transform. Thus for 3-folds $X$, 
it remains to prove Conjecture \ref{D:1.2} when $q(X)=0$. 
On the other hand, by the same technique, we can show that 
Conjecture \ref{D:1.2} holds for such varieties which are 
birational to an abelian variety (cf.\ \cite{xiep}).

(3) Assume further that $B=0$ and that $X$ is a terminal 3-fold. 
Then Conjecture \ref{D:1.2} holds provided that the second Chern 
class $c_2(X)$ is pseudo-effective (cf.\ \cite[Proposition 4.3]{xie}).
\end{rem}

In this paper, we shall consider Problem \ref{D:1.1} 
for algebraic surfaces in positive characteristic. 
There are some motivations to deal with this case. 
Firstly, both the Kodaira type vanishing theorems 
and the semipositivity theorem do not hold in general. 
Secondly, as for index 1 cover, the same as what is true in 
$\ch(k)=0$ can be false in $\ch(k)>0$. 
For instance, locally, Kawamata gave counterexamples which show 
that the index 1 cover of a log terminal surface is not necessarily 
of canonical singularities when $\ch(k)=2$ or 3 (cf.\ \cite{ka99}). 
Globally, for the Kawamata-Viehweg vanishing, the logarithmic 
Koll\'ar vanishing (see below), and the logarithmic semipositivity, 
there are counterexamples on ruled surfaces (cf.\ Examples 
\ref{D:3.6},\ref{D:3.7},\ref{D:3.9},\ref{D:3.10}). 
Thirdly, there are several kinds of pathological surfaces 
appearing in the classification theory.

We recall the logarithmic Koll\'ar vanishing theorem 
for the convenience of the reader (cf.\ \cite[Theorem 10.19]{ko95}).

\begin{thm}[Logarithmic Koll\'ar Vanishing]\label{D:1.12}
Let $f:X\rightarrow Y$ be a surjective morphism between normal proper 
varieties over an algebraically closed field $k$ with $\ch(k)=0$. 
Let $B$ be an effective $\Q$-divisor on $X$ such that $(X,B)$ is KLT. 
Let $D$ be a Cartier divisor on $X$, and $M$ a nef and big $\Q$-Cartier 
$\Q$-divisor on $Y$, such that $D\equiv K_X+B+f^*M$. Then 
$H^i(Y,R^jf_*\OO_X(D))=0$ for any $i>0$ and any $j\geq 0$.
\end{thm}

The following are the main theorems in this paper, which give a partial 
answer to Problem \ref{D:1.1} for algebraic surfaces in positive 
characteristic.

\begin{thm}\label{D:main1}
With the same assumptions as in Problem \ref{D:1.1}, 
assume further that $\dim X=2$ and $\ch(k)>0$. Then we have

(1) $H^0(X,D)\neq 0$ holds except possibly in the following cases:

\hskip 6mm {\rm (C)} $X$ is a ruled surface with $h^1(\OO_X)\geq 2$;

\hskip 6mm {\rm (D-I)} $X$ is a quasi-elliptic surface with $\chi(\OO_X)<0$;

\hskip 6mm {\rm (D-II)} $X$ is a surface of general type with $\chi(\OO_X)<0$.

(2) In Case {\rm (C)}, $H^0(X,2D)\neq 0$ always holds. 
Furthermore, if either $X$ is relatively minimal or $D$ is not big, 
then $H^0(X,D)\neq 0$ holds.
\end{thm}

\begin{thm}\label{D:main2}
For the Kawamata-Viehweg vanishing, the logarithmic Koll\'ar vanishing 
and the logarithmic semipositivity, there are counterexamples 
on ruled surfaces in any positive characteristic.
\end{thm}

We always work over an algebraically closed field $k$ of 
characteristic $p>0$ unless otherwise stated. 
For the classification theory of surfaces in positive 
characteristic, we refer the reader to \cite{mu69,bm} or \cite{ba}. 
For the definitions and results related to the minimal model theory, 
we refer the reader to \cite{kmm,km}. We use $\equiv$ to denote 
numerical equivalence, $\sim_\Q$ to denote $\Q$-linear equivalence, 
and $[B]=\sum [b_i] B_i$ to denote the round-down of a $\Q$-divisor 
$B=\sum b_iB_i$.

\begin{ack}
I would like to express my gratitude to Professors Yujiro Kawamata 
and Takao Fujita for their valuable advices and warm encouragements. 
I would also like to thank Professors Keiji Oguiso and Natsuo Saito 
for stimulating discussions. I am very grateful to the referee 
for his useful suggestions and comments. 
This work was partially supported by JSPS grant no.\ P05044.
\end{ack}

\section{Reduction to Cases}\label{D:S2}

First of all, we give an easy reduction to Problem \ref{D:1.1} 
in the surface case.

\begin{prop}\label{D:2.1}
With the same assumptions as in Problem \ref{D:1.1}, assume that 
$\dim X=2$. Then we may assume that $X$ is smooth projective, 
$B$ is a $\Q$-divisor and $H$ is ample.
\end{prop}

\begin{proof}
Let $f:Y\rightarrow X$ be the minimal resolution of $X$. 
We may write $K_Y=f^*K_X+\sum a_iE_i$, where $E_i$ are exceptional 
curves of $f$ and $-1<a_i\leq 0$ for all $i$. 
Let $B'=f^*B-\sum a_iE_i\geq 0$. Then $K_Y+B'=f^*(K_X+B)$. 
It is easy to see that $(Y,B')$ is also KLT. 
Note that $H'=f^*D-(K_Y+B')$ is nef and big, and $H^0(X,D)\neq 0$ 
is equivalent to $H^0(Y,f^*D)\neq 0$. On the other hand, 
by Kodaira's Lemma, we may assume that $B'$ is a $\Q$-divisor and 
$H'$ is ample by adding a sufficiently small $\R$-divisor to $B'$.
\end{proof}

Therefore we consider the following problem in what follows.

\begin{prob}\label{D:2.2}
Let $X$ be a smooth projective surface over an algebraically 
closed field $k$ of characteristic $p>0$, 
$B=\sum_{i=1}^m b_iB_i$ an effective $\Q$-divisor on $X$ 
such that $(X,B)$ is KLT. Let $D$ be a nef divisor on $X$ 
such that $H=D-(K_X+B)$ is ample. Does $H^0(X,D)\neq 0$ hold?
\end{prob}

Secondly, we have the following easy criterion for non-vanishing.

\begin{lem}\label{D:2.3}
If $\chi(X,D)>0$, then $h^0(X,D)>0$.
\end{lem}

\begin{proof}
We have that $h^2(X,D)=h^0(X,K_X-D)=h^0(X,-H-B)=0$ by Serre duality, 
hence the conclusion is obvious.
\end{proof}

\begin{case}[A]
$D\equiv 0$, hence $-(K_X+B)$ is ample.
\end{case}

It follows from Serre duality that $h^2(X,\OO_X)=h^0(X,K_X)=0$. We 
shall show that $h^1(X,\OO_X)=0$ by the following two lemmas.

\begin{lem}\label{D:2.4}
Assume that we are in Case (A). Then $\overline{NE}(X)=\sum\R_+[l_i]$, 
where $l_i$ are rational curves on $X$ (not necessarily extremal).
\end{lem}

\begin{proof}
For any $C\not\subset\Supp B$, we have 
$-K_X.C>B.C\geq 0$. On the other hand, we have $-K_X.B_i>
(\sum b_jB_j).B_i\geq b_iB_i^2$. If $B_i^2\geq 0$, then 
$-K_X.B_i>0$. If $B_i^2<0$, then $-(K_X+b_iB_i).B_i>
(\sum_{j\neq i}b_jB_j).B_i\geq 0$ and $2-2p_a(B_i)=-(K_X+B_i).B_i>
-(K_X+b_iB_i).B_i>0$. Hence $p_a(B_i)=0$ and $B_i\cong\PP^1$.

By permutation of the indices, we may assume that 
$B_i^2<0$ for $1\leq i\leq s$, where $0\leq s\leq m$. 
By the cone theorem (cf.\ \cite[Theorem 1.4]{mo}), we have
\[ \overline{NE}(X)=\sum_{i=1}^r\R_+[l_i]+
\overline{NE}_{K_X+\varepsilon L\geq 0}(X), \]
where $l_1,\cdots,l_r$ are extremal rational rays and $L=-(K_X+B)$.

We claim that $\overline{NE}(X)=\sum_{i=1}^r\R_+[l_i]+\sum_{j=1}^s\R_+[B_j]$. 
Indeed, for any curve $C$, we may write $C=\lim(\sum a_il_i+\sum c_kz_k)$, 
where $a_i\geq 0, c_k\geq 0$, $z_k\in NE_{K_X+\varepsilon L\geq 0}(X)$ 
are irreducible curves on $X$, and $\lim$ means the limit of vectors 
under the usual topology of $\overline{NE}(X)$. 
By definition, for each $k$ we have
\begin{eqnarray}
(K_X-\varepsilon(K_X+B)).z_k & = & (1-\varepsilon)K_X.z_k-
\varepsilon B.z_k\geq 0, \nonumber \\
K_X.z_k & \geq & \frac{\varepsilon}{1-\varepsilon}B.z_k. \nonumber
\end{eqnarray}
If $z_k\not\subset\Supp B$, then $K_X.z_k\geq 0$, 
a contradiction. Hence $z_k=B_j$ for some $1\leq j\leq n$. 
If $B_j^2\geq 0$, then $K_X.z_k\geq 0$, a contradiction. 
Hence $z_k=B_{j}\cong\PP^1$ for some $1\leq j\leq s$. 
Therefore 
$C=\lim(\sum_{i=1}^r a_il_i+\sum_{j=1}^s c_jB_j)=
\sum_{i=1}^r \bar{a}_il_i+\sum_{j=1}^s \bar{c}_jB_j$.
\end{proof}

\begin{lem}\label{D:2.5}
Assume that we are in Case (A). Let $\alpha: X\rightarrow A$ 
be the Albanese map of $X$. Then $q(X):=\dim A=0$ and 
$h^1(\OO_X)=0$.
\end{lem}

\begin{proof}
Let $M$ be an ample divisor on $A$. 
By Lemma \ref{D:2.4}, for any curve $C$ on $X$, 
we may write $C\equiv\sum a_il_i$, where $a_i\geq 0$ and $l_i$ are rational 
curves on $X$. Since $A$ contains no rational curves, $\alpha(l_i)$ 
is a point for each $i$. Then 
\[ \alpha^*M.C=\alpha^*M.(\sum a_il_i)=\sum a_i\alpha^*M.l_i=0, \]
hence $\alpha(C)$ is also a point. Thus $\alpha$ is constant and $q(X)=0$.

Note that the following inequalities hold (cf.\ \cite{bm}):
\[ 0\leq h^1(\OO_X)-q(X)\leq p_g(X)=h^2(\OO_X)=0. \]
Hence $h^1(\OO_X)=q(X)=0$.
\end{proof}

In total, in Case (A), we have $\chi(X,D)=\chi(\OO_X)=1>0$. 
As a corollary, we know that any smooth projective surface with a 
log Fano structure is rational.

\begin{case}[B]
$D\not\equiv 0$ and either

{\rm (I)} $\kappa(X)\geq 0$ and $\chi(\OO_X)\geq 0$, or

{\rm (II)} $X$ is a ruled surface and $q(X)\leq 1$, hence $\chi(\OO_X)\geq 0$.
\end{case}

In Case (B), by the Riemann-Roch theorem, we have
\begin{eqnarray}
\chi(X,D) & = & \frac{1}{2}D(D-K_X)+\chi(\OO_X) \nonumber \\
          & = & \frac{1}{2}D(H+B)+\chi(\OO_X)>0. \nonumber
\end{eqnarray}

Let us consider the remaining cases. Assume that $X$ is not contained 
in Cases (A) or (B). Let $Y$ be a relatively minimal model of $X$. 
If $\kappa(Y)=-\infty$, then $Y$ must be a $\PP^1$-bundle with $c_2(Y)<0$,
which is Case (C).

\begin{case}[C]
$D\not\equiv 0$. There exist a smooth curve $C$ with $g(C)\geq 2$ 
and a surjective morphism $f:X\rightarrow C$ such that $X$ is a 
ruled surface over $C$.
\end{case}

In characteristic zero, it is well-known that if $c_2(X)<0$, 
then $X$ is ruled. A similar result holds in positive characteristic 
due to Raynaud and Shepherd-Barron (cf. \cite[Theorem 7]{sb}).

\begin{thm}\label{D:2.6}
Let $X$ be a smooth surface over an algebraically closed field $k$ 
of positive characteristic. If $c_2(X)<0$, then $X$ is uniruled. 
In fact, there exist a smooth curve $C$ and a surjective morphism 
$f:X\rightarrow C$ such that the geometric generic fiber 
of $f$ is a rational curve.
\end{thm}

If $\kappa(Y)=0$, then $c_2(Y)\geq 0$ by the explicit classification 
(cf.\ \cite{bm}), hence $\chi(\OO_X)=\chi(\OO_Y)\geq 0$, such $X$ are 
contained in Case (B-I).
If $\kappa(Y)=1$ and $c_2(Y)=12\chi(\OO_Y)<0$, then $Y$ must be a 
quasi-elliptic surface by the classification theory and Theorem 
\ref{D:2.6}. The last one is the case that $X$ is of general type with 
$\chi(\OO_X)<0$. Therefore we have the following Case (D).

\begin{case}[D]
$D\not\equiv 0$. There exist a smooth curve $C$ and a surjective 
morphism $f:X\rightarrow C$ such that $\chi(\OO_X)<0$ and either 

{\rm (I)} the geometric generic fiber of $f$ is a rational curve 
with an ordinary cusp, or 

{\rm (II)} the geometric generic fiber of $f$ is a rational curve, 
and $X$ is of general type.
\end{case}

In $\ch(k)=0$, Case (D) cannot occur, and Case (C) is settled by 
Kawamata by using the logarithmic semipositivity theorem 
(cf.\ \cite[Theorem 3.1]{ka00}). 
Note that Case (D-I) can occur only if $\ch(k)=2$ or 3 
(cf.\ \cite{bm}), and the explicit examples have been given by 
Raynaud and Lang (cf.\ \cite{ra,la}). For Case (D-II), we can restrict 
our attention to a small class by \cite[Theorem 8]{sb}, however 
no example is known so far.

We shall discuss Case (C) in \S\ref{D:S3} and \S\ref{D:S4}.

\section{Some Counterexamples}\label{D:S3}

When $\ch(k)=p>0$, it is well-known that the Kodaira vanishing does 
not hold on surfaces in general. However, the Kodaira vanishing does 
hold on ruled surfaces, which was first proved by Tango 
(cf.\ \cite{ta72b}). In fact, we have the following theorem 
given by Mukai (cf.\ \cite{mu79}).

\begin{thm}\label{D:3.1}
Let $X$ be a smooth projective surface over an algebraically closed 
field $k$ of characteristic $p>0$. If the Kodaira vanishing does not 
hold on $X$, then $X$ must be a quasi-elliptic surface or a surface 
of general type.
\end{thm}

Furthermore, we may ask whether the Kawamata-Viehweg vanishing 
holds on ruled surfaces. This problem is important because the 
Kawamata-Viehweg vanishing gives a sufficient condition for 
the effective non-vanishing in Case (C). Roughly speaking, 
the vanishing of $H^1(X,D)$ implies the non-vanishing of 
$H^0(X,D)$ by virtue of the Fourier-Mukai transform. 
This idea was first used in \cite{ch,ccz}. 
We recall the following theorem due to Mukai 
(cf.\ \cite[Theorem 2.2]{mu81}).

\begin{thm}\label{D:3.2}
Let $A$ be an abelian variety, $\hat A$ its dual abelian variety, 
$\cP$ the Poincar\'e line bundle on $A\times\hat A$. 
Then the Fourier-Mukai transform $\Phi^{\cP}_{A\to\hat A}:
D(A)\rightarrow D(\hat A)$, $\FF^\bullet\mapsto {\mathbf R}\pi_{{\hat A}^*}
(\pi_A^*\FF^\bullet\stackrel{\mathbf L}{\otimes}\cP)$ is an equivalence of 
derived categories.

Let $\FF$ be a coherent sheaf on $A$. Assume that 
$H^i(A,\FF\otimes P)=0$ for all $P\in\Pic^0(A)$ and all $i\neq i_0$. 
Then the dual sheaf $\hat\FF=\Phi^{\cP}_{A\to\hat A}(\FF)$ 
is a locally free sheaf on $\hat A$ of rank $h^{i_0}(A,\FF)$.
\end{thm}

\begin{prop}\label{D:3.3}
Assume that we are in Case (C), and that $H^1(X,D+f^*P)=0$ 
for any $P\in\Pic^0(C)$. Then $H^0(X,D)\neq 0$ holds.
\end{prop}

\begin{proof}
Let $\alpha:X\rightarrow A=\Alb(X)$ be the Albanese map of $X$. 
Then $\alpha(X)=C\subset A$. Let $\FF=\alpha_*\OO_X(D)$ be the 
coherent sheaf on $A$. Then we have that $H^i(X,D+\alpha^*P)=
R^i\alpha_*(D+\alpha^*P)=0$ for any $P\in\Pic^0(A)$ and any $i>0$ 
by the assumption and easy computations. 
It follows from the Leray spectral sequence that 
$H^i(A,\FF\otimes P)=0$ for any $P\in\Pic^0(A)$ and any $i>0$, 
hence by Theorem \ref{D:3.2}, 
its dual $\hat\FF$ is a locally free sheaf of rank 
$h^0(A,\FF)=h^0(X,D)$. If $H^0(X,D)=0$, then $\hat\FF=0$, 
hence $\FF=0$. Next we prove that $\FF\neq 0$. Consider the 
general fiber $F$ of $f:X\rightarrow C$, then the stalk of $\FF$ 
at the general point of $C$ is isomorphic to $H^0(F,D|_F)\neq 0$ 
since $D$ is nef and $F\cong\PP^1$.
\end{proof}

\begin{rem}\label{D:3.4}
Proposition \ref{D:3.3} gives a new proof of the surface case 
of Conjecture \ref{D:1.2}, since the Kawamata-Viehweg vanishing 
theorem holds in characteristic zero.
\end{rem}

Even if the Kodaira vanishing holds on ruled surfaces, 
we cannot expect that the Kawamata-Viehweg vanishing holds 
on ruled surfaces in general. Next we shall give some 
counterexamples for the Kawamata-Viehweg vanishing 
on ruled surfaces. The constructions are similar to, however 
generalize those to some extent, which were given by Raynaud 
to yield the counterexamples for the Kodaira vanishing on 
quasi-elliptic surfaces and general type surfaces (cf.\ \cite{ra}).

\begin{defn}\label{D:3.5}
Let $C$ be a smooth projective curve over an algebraically 
closed field $k$ of characteristic $p>0$. Let $f\in K(C)$ 
be a rational function on $C$. 
\[ n(f):=\deg\biggl{[}\frac{(df)}{p}\biggr{]}, \]
where $(df)=\sum_{x\in C}v_x(df)x$ is the divisor associated to the 
rational differential 1-form $df$. 
We denote $K^p(C)=\{ f^p \, | \, f\in K(C) \}$.
\[ n(C):=\max\{n(f) \, | \, f\in K(C),f\not\in K^p(C) \}. \]

If $f\not\in K^p(C)$, then $(df)$ is a canonical divisor on $C$ 
with degree $2(g-1)$. It is easy to see that $n(C)\leq [2(g-1)/p]$.
\end{defn}

\begin{ex}\label{D:3.6}
There do exist smooth projective curves $C$ such that $n(C)>0$ 
for each characteristic $p>0$.

(1) Let $h\geq 3$ be an odd integer, $p\geq 3$. Let $C$ be the 
projective completion at infinity of the affine curve defined by 
$y^2=x^{ph}+x^{p+1}+1$. It is easy to verify that $C$ is a 
smooth hyperelliptic curve and that $(d(y/x^p))=(ph-3)z_\infty$, 
where $z_\infty$ is the infinity point of $C$ 
(cf.\ \cite[Ch.\ III, \S 6.5]{sha}). Hence $n(C)=n(y/x^p)=h-1>0$.

(2) (cf.\ \cite{ra}). 
Let $h>2$ be an integer. Let $C$ be the projective completion 
at infinity of the Artin-Schreier cover of the affine line 
defined by $y^{hp-1}=x^p-x$. It is easy to verify that $C$ 
is a smooth curve of genus $g$ with $2(g-1)=p(h(p-1)-2)$, 
and that $(dy)=p(h(p-1)-2)z_\infty$, where $z_\infty$ is 
the infinity point of $C$. Hence $n(C)=n(y)=h(p-1)-2>0$.

(3) (cf.\ \cite{ta72a}). Let $C\subset\PP^2$ be the curve defined 
by $x_0^{p+1}=x_1x_2(x_0^{p-1}+x_1^{p-1}-x_2^{p-1})$, where $p\geq 3$. 
We can show that $C$ is smooth and that $n(C)=n(x_0/x_1)=p-2>0$.
\end{ex}

\begin{ex}\label{D:3.7}
{\it Let $C$ be a smooth projective curve over an algebraically 
closed field $k$ of characteristic $p>0$. If $n(C)>0$, 
then there are a $\PP^1$-bundle $f:X\rightarrow C$, 
an effective $\Q$-divisor $B$ and an integral divisor $D$ 
on $X$ such that $(X,B)$ is KLT and $H=D-(K_X+B)$ is ample. 
However $H^1(X,D)\neq 0$. }

Let $F:C\rightarrow C$ be the Frobenius map. We have the following 
exact sequences of $\OO_C$-modules:
\begin{eqnarray}
0\rightarrow \OO_C\rightarrow F_*\OO_C\rightarrow \BB^1\rightarrow 0 
\label{es:1} \\
0\rightarrow \BB^1\rightarrow F_*\Omega^1_C\stackrel{c}{\rightarrow} 
\Omega^1_C\rightarrow 0 \label{es:2} 
\end{eqnarray}
where $\BB^1$ is the image of the map $F_*(d):F_*\OO_C\rightarrow 
F_*\Omega^1_C$, and $c$ is the Cartier operator (cf.\ \cite{ta72a}).

Let $\LL=\OO_C(L)$ be a line bundle on $C$. 
Tensor (\ref{es:2}) by $\OO_C(-L)$, we have:
\[ 0\rightarrow \BB^1(-L)\rightarrow F_*(\Omega^1_C(-F^*L))
\stackrel{c(-L)}{\rightarrow} \Omega^1_C(-L)\rightarrow 0. \]
Thus $H^0(C,\BB^1(-L))=\bigl\{ df \, | \, f\in K(C), \, 
(df)\geq pL \big\}$. 
Since $n(C)>0$, there exists an $f_0\in K(C)$ such that 
$n(f_0)=\deg [(df_0)/p]=n(C)>0$. Let $L=[(df_0)/p]$. Then $\deg L=n(C)>0$ 
and $(df_0)\geq pL$, hence $0\neq df_0\in H^0(C,\BB^1(-L))$, 
and we can regrad the line bundle $\LL=\OO_C(L)\subset\BB^1$.

Tensor (\ref{es:1}) by $\LL^{-1}$ and take cohomology, we have:
\[ 0\rightarrow H^0(C,\BB^1(-L))\stackrel{\eta}{\rightarrow} 
H^1(C,\LL^{-1})\stackrel{F^*}{\rightarrow} H^1(C,\LL^{-p}). \]
Since $\eta$ is injective, we may take the element $0\neq\eta(df_0)
\in H^1(C,\LL^{-1})$, which determines the following extension sequence: 
\begin{eqnarray}
0\rightarrow \OO_C\rightarrow \EE\rightarrow \LL\rightarrow 0. \label{es:3}
\end{eqnarray}
Pull back the exact sequence (\ref{es:3}) by the Frobenius map $F$, 
we have the following split exact sequence:
\begin{eqnarray}
0\rightarrow \OO_C\rightarrow F^*\EE\rightarrow \LL^p\rightarrow 0, 
\label{es:4}
\end{eqnarray}
since the obstruction of extension of (\ref{es:4}) is just 
$F^*\eta(df_0)=0$.

Let $X=\PP(\EE)$ be the $\PP^1$-bundle over $C$, $f:X\rightarrow C$ 
the projection, $\OO_X(1)$ the tautological line bundle. The sequence 
(\ref{es:3}) determines a section $E$ of $f$ such that 
$\OO_X(E)\cong \OO_X(1)$, and $E$ corresponds to a section 
$s\in H^0(X,\OO_X(1))=H^0(C,\EE)$ which is the image of 1 under the map 
$H^0(C,\OO_C)\hookrightarrow H^0(C,\EE)$. The sequence (\ref{es:4}) 
induces an exact sequence:
\[ 0\rightarrow \OO_C\rightarrow F^*\EE\otimes\LL^{-p}\rightarrow 
\LL^{-p}\rightarrow 0, \]
which determines a section $t\in H^0(X,\OO_X(p)\otimes f^*\LL^{-p})$ 
through the maps $H^0(C,\OO_C)\hookrightarrow H^0(C,F^*\EE\otimes\LL^{-p})
\hookrightarrow H^0(C,S^p(\EE)\otimes\LL^{-p})=
H^0(X,\OO_X(p)\otimes f^*\LL^{-p})$. The section $t$ determines an 
irreducible curve $C'$ on $X$ such that $\OO_X(C')\cong 
\OO_X(p)\otimes f^*\LL^{-p}$. It is easy to verify that 
both $E$ and $C'$ are smooth over $k$, and $E\cap C'=\emptyset$.

(\dag) {\it Assume that $p\geq 3$}.

Let $B=\frac{1}{2}C'$, $D=K_X+\frac{p+1}{2}E+\frac{1-p}{2}f^*L=
\frac{p-3}{2}E+f^*(K_C+\frac{3-p}{2}L)$. Then 
$H=D-(K_X+B)=\frac{1}{2}(E+f^*L)$. 
It is easy to see that $(X,B)$ is KLT. Since $E^2=\deg\EE=\deg\LL>0$, 
$E$ is a nef divisor on $X$. On the other hand, $E$ is $f$-ample, hence 
$H$ is an ample $\Q$-divisor on $X$. Next we show that $H^1(X,D)\neq 0$.

Consider the Leray spectral sequence $E_2^{i,j}=H^i(C,R^jf_* ?)
\Rightarrow H^{i+j}(X,?)$. Since $E_2^{i,j}=0$ for $i\geq 2$, 
by the five term exact sequence we have
\[ H^1(X,D)\cong H^1(X,-H-B)^\vee\supseteq 
H^0(C,R^1f_*\OO_X(-\frac{p+1}{2})\otimes \LL^{\frac{p-1}{2}})^\vee \]
By the relative Serre duality, 
\begin{eqnarray}
R^1f_*\OO_X(-\frac{p+1}{2})^\vee & \cong &
f_*(\OO_X(\frac{p+1}{2})\otimes\omega_{X/C}) \nonumber \\
& = & f_*(\OO_X(\frac{p-3}{2})\otimes f^*\LL)=
S^{(p-3)/2}(\EE)\otimes\LL. \nonumber
\end{eqnarray}
Since $S^{(p-3)/2}(\EE)$ has a quotient sheaf $\LL^{(p-3)/2}$, 
$R^1f_*\OO_X(-(p+1)/2)$ has a subsheaf $\LL^{-(p-1)/2}$. 
Thus $H^1(X,D)\supseteq H^0(C,\OO_C)^\vee=k$, which is desired.

(\ddag) {\it Assume that $p=2$}.

Let $B=\frac{2}{3}C'$, $D=K_X+2E-f^*L=f^*K_C$, 
$H=D-(K_X+B)=\frac{1}{3}(2E+f^*L)$. 
It is easy to verify that $(X,B)$ is KLT and $H$ is $\Q$-ample. 
By the same argument, we have 
$H^1(X,D)=H^1(X,f^*\omega_C)\cong H^1(X,\omega_{X/C})^\vee
\supseteq H^0(C,R^1f_*\omega_{X/C})^\vee\cong H^0(C,\OO_C)^\vee=k$, 
which is desired. \qed
\end{ex}

By Examples \ref{D:3.6}-7, there do exist counterexamples 
for the Kawamata-Viehweg vanishing on ruled surfaces. Furthermore, 
it is easy to verify that $D$ is nef and $|D|\neq\emptyset$ 
in both cases. Hence it follows that the Kawamata-Viehweg vanishing 
is a sufficient but not a necessary condition for the effective 
non-vanishing in Case (C).

Examples \ref{D:3.6}-7 also give the counterexamples for the 
$\Q$-divisor version of the Kawamata-Viehweg vanishing 
(cf.\ \cite[Theorem 1-2-3]{kmm}). Indeed, we can take $D-(K_X+B)$ 
as the required $\Q$-divisor. However, it is unknown whether there 
exist counterexamples for its nef and big version mentioned below. 
So it is interesting to take the following problem into 
account, which is compared with Theorem \ref{D:3.1}.

\begin{prob}\label{D:3.8}
Let $X$ be a smooth projective surface over an algebraically 
closed field $k$ of characteristic $p>0$. Let $D$ be an integral 
divisor on $X$ such that $D-K_X$ is nef and big. Assume that 
$X$ is neither quasi-elliptic nor of general type. 
Does $H^1(X,D)=0$ hold?
\end{prob}

In arbitrary characteristic, the Koll\'ar vanishing of $f_*\omega_X$ 
for ruled surfaces is trivial (cf.\ \cite[Theorem 2.1]{ko}). However, 
we shall give counterexamples for the logarithmic Koll\'ar vanishing 
in positive characteristic (cf.\ Theorem \ref{D:1.12}).

\begin{ex}\label{D:3.9}
{\it Let $C$ be a smooth projective curve over an algebraically 
closed field $k$ of characteristic $p>0$. If $n(C)>0$, 
then there are a $\PP^1$-bundle $f:X\rightarrow C$, 
an effective $\Q$-divisor $B'$ and an integral divisor $D$ 
on $X$ such that $(X,B')$ is KLT and $D\sim_\Q K_X+B'+f^*M$, 
where $M$ is an ample $\Q$-divisor on $C$. 
However $H^1(C,f_*\OO_X(D))\neq 0$. }

It is just Example \ref{D:3.7}. We use the same notation and 
assumptions as in Example \ref{D:3.7}.
When $p\geq 3$, let $B'=\frac{1}{2}(E+C')$ and $M=\frac{1}{2}L$. 
When $p=2$, let $B'=\frac{2}{3}(E+C')$ and $M=\frac{1}{3}L$. 
Then $D\sim_\Q K_X+B'+f^*M$. It follows from $R^1f_*\OO_X(D)=0$ 
and the Leray spectral sequence that 
$H^1(C,f_*\OO_X(D))\cong H^1(X,D)\neq 0$. \qed
\end{ex}

In characteristic zero, Kawamata settled Case (C) by means of 
the logarithmic semipositivity theorem (cf.\ Theorem \ref{D:1.10}). 
In arbitrary characteristic, the semipositivity of $f_*\omega_{X/C}$ 
is trivial for ruled surfaces, however we shall give counterexamples 
for the logarithmic semipositivity in positive characteristic.

\begin{ex}\label{D:3.10}
{\it Let $C$ be a smooth projective curve over an algebraically 
closed field $k$ of characteristic $p>0$. Assume that $n(C)>0$ and 
that the following condition holds:

$(*)$ If $p=2$, then $\frac{1}{3}L$ is integral, 
where $L=[(df_0)/p]$ is the divisor on $C$ for some rational 
function $f_0\in K(C)$ such that $n(f_0)=n(C)$.

Then there are a $\PP^1$-bundle $f:X\rightarrow C$, 
an effective $\Q$-divisor $B'$ and an integral divisor $D'$ on $X$ 
such that $(X,B')$ is KLT and $D'\sim_\Q K_{X/C}+B'$. 
However $f_*\OO_X(D')$ is not semipositive. }

When $p\geq 5$, the counterexample is just Example \ref{D:3.7}. 
We use the same notation and assumptions as in Example \ref{D:3.7}. 
Since $H$ is ample, we can take a general member $M\in |nH|$ for 
$n$ sufficiently large and divisible such that $M$ is irreducible and 
smooth, and $B'=B+\frac{1}{n}M$ has simple normal crossing support, 
hence $(X,B')$ is KLT. Let $D'=D-f^*K_C$. Then $D'\sim_\Q K_{X/C}+B'$. 
Since $D'|_F=D|_F$ is nef hence basepoint free on $F\cong\PP^1$, 
the canonical homomorphism $f^*f_*\OO_X(D') 
\twoheadrightarrow \OO_X(D')$ is surjective. If $f_*\OO_X(D')$ 
were semipositive, then we would have $D'=D-f^*K_C$ is nef on $X$. 
However, $(D-f^*K_C)C'=\frac{3-p}{2}f^*L.C'<0$, 
this is absurd. Hence $f_*\OO_X(D')$ never be semipositive.

For $p<5$, we need to modify Example \ref{D:3.7} slightly. 
When $p=3$, let $B=\frac{5}{6}C'$, $D=E+f^*(K_C-L)$. 
Then $(X,B)$ is KLT and $H=D-(K_X+B)=\frac{1}{2}(E+f^*L)$ is 
ample. However $D'=D-f^*K_C$ satisfies $D'.C'=(E-f^*L)C'=-f^*L.C'<0$, 
hence $f_*\OO_X(D')$ never be semipositive.

When $p=2$, we need the additional assumption $(*)$ mentioned above. 
Let $B=\frac{5}{6}C'$, $D=E+f^*(K_C-\frac{1}{3}L)$. Then 
$(X,B)$ is KLT and $H=D-(K_X+B)=\frac{1}{3}(4E+f^*L)$ is 
ample. However $D'=D-f^*K_C$ satisfies $D'.C'=
(E-\frac{1}{3}f^*L)C'=-\frac{1}{3}f^*L.C'<0$, 
hence $f_*\OO_X(D')$ never be semipositive. \qed
\end{ex}

Note that the assumption $(*)$ can be realized by 
Example \ref{D:3.6}(2). Indeed, let $h-2$ be a positive 
integer divisible by 3, then we are done. 
Hence there do exist counterexamples for the logarithmic 
semipositivity on ruled surfaces. 

Furthermore, it is easy to verify that $D$ is nef and 
$|D|\neq\emptyset$ in both cases, hence it follows that 
the logarithmic semipositivity is a sufficient but not a necessary 
condition for the effective non-vanishing in Case (C).

Let us compare the two approaches for proving the effective 
non-vanishing conjecture for surfaces in characteristic zero. 
Of course, we only need to treat Case (C). 
Since the semipositivity theorem can be deduced from the 
Koll\'{a}r vanishing theorem (cf.\ \cite[Corollary 3.7]{ko}), 
the approach provided by Kawamata gives the diagram (1), and 
Proposition \ref{D:3.3} gives the diagram (2) as follows:
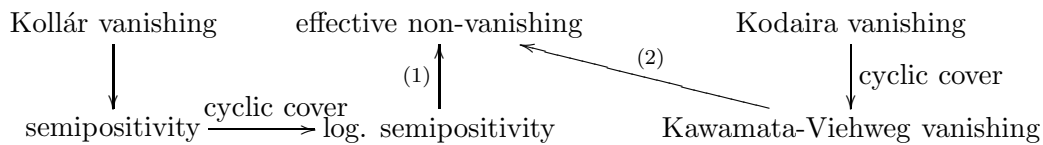
\[ \xymatrix{
\kolv \ar[d] & \en & \kodv \ar[d]^{\cc} \\
\spa \ar[r]^{\cc} & \lsp \ar[u]^{(1)} & \kvv \ar[ul]_{(2)}
} \]

In characteristic zero, the vanishing theorem is the start point 
of both approaches, and the cyclic cover trick plays a more important 
role in both proofs. However, Examples \ref{D:3.7} and \ref{D:3.10} 
show that, to some extent, the cyclic cover trick does not behave well 
in positive characteristic. It will turn out in the next section 
that without the cyclic cover trick, we could not deal with the 
case $B\neq 0$ effectively.

\section{Ruled Surface Case}\label{D:S4}

Firstly, there is a partial answer to the effective 
non-vanishing in Case (C), whose proof is numerical, hence valid 
in positive characteristic (cf.\ \cite[Proposition 4.1(2a)]{am}).

\begin{prop}\label{D:4.1}
Let $F$ be the general fiber of $f:X\rightarrow C$. 
If $H.F>1$, then $H^0(X,D)\neq 0$ (This is true even if 
$H=D-(K_X+B)$ is nef and big).
\end{prop}

Proposition \ref{D:4.1} guarantees the non-vanishing for the absolute 
case, i.e.\ $B=0$ and $D-K_X$ is nef and big, since $H.F\geq -K_X.F=2>1$. 
Hence we have to consider the case $B\neq 0$.

Secondly, we shall prove the following theorem as a first step.

\begin{thm}\label{D:4.2}
In Case (C), assume furthermore that $X$ is relatively minimal. 
Then $H^0(X,D)\neq 0$.
\end{thm}

Let us fix some notation. Assume that $f:X=\PP(\EE)\rightarrow C$ 
is a $\PP^1$-bundle over $C$ associated to a normalized rank 2 
locally free sheaf $\EE$ on $C$. Let $e=-\deg\EE$, 
$E$ the canonical section of $f$ with $E^2=-e$, $F$ the fiber of $f$. 
Note that the proof of Theorem \ref{D:4.2} is also numerical, and that 
we only need the condition $[B]=0$, so the KLT assumption of $(X,B)$ 
is unnecessary.

Assume that $e\geq 0$. It is easy to see that if $L\equiv aE+bF$ 
is an irreducible curve on $X$, then either $L=E,F$ 
or $a>0,b\geq ae\geq 0$. Hence $L^2=a(2b-ae)\geq 0$ in the latter case. 
In other words, if $L^2<0$ then $L=E$ and $e>0$. We may write 
$B=aE+B'$ with $E\not\subset \Supp B'$. Then $B'$ is nef, 
$H+B'=D-(K_X+aE)$ is ample and $(H+B').F=(D-K_X-aE).F\geq 2-a>1$. 
By Proposition \ref{D:4.1}, we have $H^0(X,D)\neq 0$.

It remains to deal with the case $e<0$. Let $B=\sum_{i\in I}b_iB_i$. 
If $B_i^2\geq 0$, then $B_i$ is a nef divisor on $X$, and we can 
move $b_iB_i$ from $B$, add $b_iB_i$ to $H$ and keep $D$ unchanged 
to consider the non-vanishing problem. Hence we may assume that 
$B_i^2<0$ for all $i\in I$. Since $B_i$ are numerically independent 
and $\rho(X)=2$, we have $|I|\leq 1$. Indeed, if $B_1,B_2$ are 
distinct components of $B$, then we may write $F\equiv c_1B_1+c_2B_2$, 
where $c_i$ are rational numbers and at least one of $c_i$ is positive. 
If both $c_i>0$, then both $B_i.F=0$, hence $B_i=F$, 
a contradiction. If $c_1>0,c_2\leq 0$, then $F.B_1=
c_1B_1^2+c_2B_2.B_1<0$, a contradiction.

Therefore we have only to consider the following case:

\begin{case}[C-M]
Let $f:X\rightarrow C$ is a $\PP^1$-bundle over a smooth curve $C$ 
of genus $g\geq 2$ with invariant $e<0$. Let $D\not\equiv 0$ be a 
nef divisor on $X$, $B=cG$, where $0<c<1$ and $G$ is an irreducible 
curve on $X$ with $G^2<0$, such that $H=D-(K_X+B)$ is ample.
\end{case}

We need an easy lemma (cf.\ \cite[Ch.\ V, Ex.\ 2.14]{ha}):

\begin{lem}\label{D:4.3}
With the same assumptions as in Case (C-M).

{\rm (i)} If $G\equiv xE+yF$ is an irreducible curve $\neq E,F$, 
then either $x=1, y\geq 0$, or $2\leq x\leq p-1, y\geq xe/2$, or 
$x\geq p, y\geq xe/2+1-g$.

{\rm (ii)} If $D\equiv aE+bF$ is ample, then $a>0, b>ae/2$.
\end{lem}

\begin{lem}\label{D:4.4}
Assume that we are in Case (C-M). Then $H^0(X,D)\neq 0$.
\end{lem}

\begin{proof}
Let $D\equiv aE+bF$, $G\equiv xE+yF$. Then $a\geq 0$, $b\geq ae/2$ 
and $x,y$ satisfy the condition mentioned in Lemma \ref{D:4.3}(i). 
We have $H\equiv aE+bF+2E+(2-2g+e)F-cxE-cyF=(a+2-cx)E+(b+2-2g+e-cy)F$. 
Since $H$ is ample, the following conditions hold by Lemma \ref{D:4.3}(ii):
\[ a+2-cx > 0, \,\, b+2-2g+e-cy > \frac{1}{2}(a+2-cx)e. \]
By Lemma \ref{D:4.3}(i) and the later inequality, we have 
\[ b-\frac{1}{2}ae>2g-2+c(y-\frac{1}{2}xe)>(2-c)(g-1)>g-1. \]
It follows from the Riemann-Roch theorem that
\[ \chi(X,D)=\frac{1}{2}D(D-K_X)+\chi(\OO_X)=(a+1)(b-\frac{1}{2}ae+1-g)>0, \]
which also completes the proof of Theorem \ref{D:4.2}.
\end{proof}

Let $X$ be a smooth projective surface, 
$B$ an effective $\Q$-divisor such that $(X,B)$ is KLT. 
Let $D$ be a nef divisor on $X$ such that 
$H=D-(K_X+B)$ is ample. Next, we consider the reduction of 
the effective non-vanishing problem for the triple $(X,B;D)$ 
under the $(-1)$-curve contractions.

Let $g:X\rightarrow Y$ be a contraction of a $(-1)$-curve $l\subset X$. 
Assume that there exists a divisor $D_Y$ on $Y$ such that $D=g^*D_Y$ 
(this condition is equivalent to $D.l=0$). It is easy to verify that 
$D_Y$ is nef. Let $B_Y:=g_*B$ be the strict transform of $B$. Then $B_Y$ 
is also an effective divisor with $[B_Y]=0$. We may write
\[ K_X+B=g^*(K_Y+B_Y)+dl, \]
where $d>0$ since $D-(K_X+B)$ is ample. It follows from $d>0$ that 
$(Y,B_Y)$ is again KLT. Let $C$ be an irreducible curve on $Y$, 
it is easy to verify that $(D_Y-(K_Y+B_Y)).C=(D-(K_X+B)).g^*C>0$ 
and that $(D_Y-(K_Y+B_Y))^2=(D-(K_X+B))^2+d^2>0$, hence 
$H_Y=D_Y-(K_Y+B_Y)$ is also ample by the Nakai-Moishezon criterion.

\begin{defn}\label{D:4.5}
Given a triple $(X,B;D)$. Let $g: X\rightarrow Y$ be a birational 
morphism to a smooth projective surface $Y$. Assume that $D$ is 
$g$-trivial, i.e.\ there exists a divisor $D_Y$ on $Y$ such that 
$D=g^*D_Y$. Then the induced triple $(Y,B_Y;D_Y)$ is called the 
reduction model of $(X,B;D)$.
\end{defn}

Note that $H^0(X,D)=H^0(Y,D_Y)$, hence the reduction model does 
give a reduction to the effective non-vanishing problem. As an 
application of Theorem \ref{D:4.2}, we know that in Case (C) 
if $(X,B;D)$ admits a relatively minimal reduction model, 
then the effective non-vanishing holds.

\begin{rem}\label{D:4.6}
In general, given a birational morphism $g: X\rightarrow Y$, 
e.g.\ $Y$ is a relatively minimal model of $X$, even if $D$ is not 
$g$-trivial, we also can define $D_Y=g_*D$ as the push-out 
of algebraic cycles. It is easy to verify that $D_Y$ is nef and 
$D_Y-(K_Y+B_Y)$ is ample. However this model is not good, 
since the pair $(Y,B_Y)$ is not necessarily KLT, and in general, 
$H^0(X,D)=H^0(Y,D_Y)$ does not hold by observing the following 
two examples.
\end{rem}

\begin{ex}\label{D:4.7}
Let $X$ be a smooth projective surface, $B_1$ a $(-1)$-curve and
$B_2,B_3$ smooth curves on $X$ such that $B_1,B_2,B_3$ 
intersect transversally at one point $p\in X$. 
Let $B=\frac{2}{5}B_1+\frac{4}{5}B_2+\frac{3}{4}B_3$. Then we can 
verify that the pair $(X,B)$ is KLT by blowing up at $p$. 
Let $g: X\rightarrow Y$ be the contraction of $B_1$. Then the pair 
$(Y,B_Y)$ is not KLT since the discrepancy of the exceptional divisor 
with center $p$ is $-\frac{11}{10}<-1$.

Let $X=\F_1$ be the Hirzebruch surface, $g:X\rightarrow Y=\PP^2$ the 
contraction of the $(-1)$-section, $D$ the fiber of $X$ over $\PP^1$. 
Then $D_Y=g_*D$ is a line in $\PP^2$. It is easy to see that 
$H^0(X,D)<H^0(Y,D_Y)$.
\end{ex}

Let us return to the argument of the effective non-vanishing problem 
on ruled surfaces.

\begin{lem}\label{D:4.8}
Assume that we are in Case (C). Let $F$ be the general fiber of 
$f: X\rightarrow C$. If $D.F\leq 1$, then $H^0(X,D)\neq 0$ holds.
\end{lem}

\begin{proof}

By Theorem \ref{D:4.2}, we may assume that $X$ is not relatively minimal.

(1) $D.F=0$

By assumption, $X$ contains a $(-1)$-curve $l$ which is contained 
in some fiber $F_0$ of $f$. The inequality $0\leq D.l\leq D.F_0=0$ 
implies that $D.l=0$. We may consider the contraction 
$g: X\rightarrow Y$ of $l$ and the reduction model $(Y,B_Y;D_Y)$. 
Since $D_Y.F=D.F=0$, finally we can obtain a relatively minimal 
reduction model by induction.

(2) $D.F=1$

By a similar way, we may contract all $(-1)$-curves $l$ with $D.l=0$, 
at last, to obtain a reduction model $(Y,B_Y;D_Y)$ such that $D_Y$ is 
positive on any $(-1)$-curve on $Y$. We claim that $Y$ is relatively 
minimal. Otherwise, there would exist a $(-1)$-curve $l_0$ contained in 
some fiber $F_0=\sum_{i=0}^r l_i$ such that all of $l_i$ are smooth 
rational curves with negative self-intersections. The inequality 
$0<D_Y.l_0\leq D_Y.F_0=D.F=1$ implies that $D_Y.l_0=1$ and $D_Y.l_i=0$ 
for all $i>0$, hence $l_i$ are not $(-1)$-curves and $K_Y.l_i\geq 0$ for 
all $i>0$. Thus we have $-2=K_Y.F_0=K_Y.l_0+\sum_{i=1}^r K_Y.l_i\geq -1$, 
a contradiction.
\end{proof}

\begin{rem}\label{D:4.9}
The proof of the case $D.F=1$ in Lemma \ref{D:4.8} has already 
appeared in that of Proposition 4.1(2b) of \cite{am}. However, 
there is a mistake in the remaining argument of the relatively 
minimal case. So we give a complete proof here for the convenience 
of the reader.
\end{rem}

Due to an idea of Ambro, we can give the following partial results.

\begin{prop}\label{D:4.10}
In Case (C), $H^0(X,2D)\neq 0$ always holds. 
Furthermore, assume that the Iitaka dimension $\kappa(X,-K_X)\geq 0$ 
or the numerical dimension $\nu(D)=1$. Then $H^0(X,D)\neq 0$ holds.
\end{prop}

\begin{proof}
By Lemma \ref{D:4.8}, we may assume that $a=D.F\geq 2$. 
Apply \cite[Lemma 4.2]{am} to $D/a$, then we have 
$\chi(\OO_X)\geq -D(D+aK_X)/2a^2$, hence
\begin{eqnarray}
\chi(X,2D) & \geq & D(2D-K_X)-\frac{1}{2a^2}D(D+aK_X) 
\nonumber \\
 & = & \frac{2a+1}{2a}D((2-\frac{1}{a})D-K_X) \nonumber \\
 & = & \frac{2a+1}{2a}D((1-\frac{1}{a})D+H+B)>0. \nonumber
\end{eqnarray}

If $\kappa(X,-K_X)\geq 0$, then we have $D.K_X\leq 0$, hence
\begin{eqnarray}
\chi(X,D) & \geq & \frac{1}{2}D(D-K_X)-\frac{1}{2a^2}D(D+aK_X) 
\nonumber \\
 & = & \frac{a+1}{2a}D((1-\frac{1}{a})D-K_X) \nonumber \\
 & = & \frac{a^2-1}{2a^2}D(H+B)-\frac{a+1}{2a^2}D.K_X>0. \nonumber
\end{eqnarray}

If $\nu(D)=1$, then $D$ is nef but not big, i.e.\ $D^2=0$. Hence 
$D.(-K_X)=D(H+B)>0$, which implies $H^0(X,D)\neq 0$.
\end{proof}

If we denote (CR) the subcase of (C) where $D$ is nef and big and 
$\kappa(X,-K_X)=-\infty$, then it remains to deal with Problem 
\ref{D:1.1} in Case (CR) and Case (D). 
It is expected that $H^0(X,D)\neq 0$ should hold in Case (CR). 
Until now, we cannot say anything for quasi-elliptic 
surfaces and general type surfaces whose Euler 
characteristics are negative. It is expected that in Case (D), 
the universal integer $m$ should be greater than 1. 
We shall treat these in a subsequent paper.

\textsc{Department of Mathematics, Tokyo Institute of Technology, 2-12-1 Oh-okayama, Meguro, Tokyo 152-8551, Japan}

\textit{E-mail address}: \texttt{xie@math.titech.ac.jp}


\small
\begin{thebibliography}{KMM87}

\bibitem[Am99]{am}
F. Ambro, Ladders on Fano varieties, 
{\it Algebraic geometry, 9. J. Math. Sci.}, 
{\bf 94} (1999), 1126--1135.

\bibitem[Ba01]{ba}
L. B${\rm \breve{a}}$descu, 
{\it Algebraic surfaces}, Springer, 2001.

\bibitem[BM]{bm}
E. Bombieri, D. Mumford, 
Enriques' classification of surfaces in char.\ p, II, 
Complex Analysis and Algebraic Geometry, A collection of papers 
dedicated to K. Kodaira, 1977, 23--42; 
Enriques' classification of surfaces in char.\ p, III, 
{\it Inv. Math.}, {\bf 36} (1976), 197--232.

\bibitem[CCZ05]{ccz}
J. A. Chen, M. Chen, D.-Q. Zhang, 
A nonvanishing theorem for $\Q$-divisors on surfaces, 
{\it J. Algebra}, {\bf 293} (2005), 363--384.

\bibitem[CH02]{ch}
J. A. Chen, C. Hacon, 
Linear series of irregular varieties,
{\it Algebraic Geometry in East Asia, Japan}, 2002, 143--153.

\bibitem[Ha77]{ha}
R. Hartshorne, {\it Algebraic geometry}, Springer-Verlag, 1977.

\bibitem[Ka99]{ka99}
Y. Kawamata, Index 1 covers of log terminal surface singularities, 
{\it J. Algebraic Geom.}, {\bf 8} (1999), 519--527.

\bibitem[Ka00]{ka00}
Y. Kawamata, On effective non-vanishing and base-point-freeness, 
{\it Asian J. Math.}, {\bf 4} (2000), 173--182.

\bibitem[KMM87]{kmm}
Y. Kawamata, K. Matsuda, K. Matsuki, 
Introduction to the minimal model problem, 
Alg. Geom. Sendai 1985, 
{\it Adv.\ Stud.\ Pure Math.}, {\bf 10} (1987), 283--360.

\bibitem[Ko86]{ko}
J. Koll\'{a}r,
Higher direct images of dualizing sheaves I, 
{\it Ann. Math.}, {\bf 123} (1986), 11--42.

\bibitem[Ko95]{ko95}
J. Koll\'{a}r,
{\it Shafarevich maps and automorphic forms}, 
Princeton Univ.\ Press, 1995.

\bibitem[KM98]{km}
J. Koll\'{a}r, S. Mori, 
{\it Birational geometry of algebraic varieties}, 
Cambridge Tracts in Math., {\bf 134} (1998).

\bibitem[La79]{la}
W. E. Lang, 
Quasi-elliptic surfaces in characteristic three, 
{\it Ann. Scient. \'Ec. Norm. Sup.}, {\bf 12} (1979), 473--500.

\bibitem[Mo82]{mo}
S. Mori, 
Threefolds whose canonical bundles are not numerically effective, 
{\it Ann. Math.}, {\bf 116} (1982), 133--176.

\bibitem[Mu69]{mu69}
D. Mumford, 
Enriques' classification of surfaces in char.\ p, I, 
Global Analysis, 1969, 325--339.

\bibitem[Mu79]{mu79}
S. Mukai, 
On counterexamples for the Kodaira vanishing theorem and the 
Yau inequality in positive characteristic (in Japanese), 
Symposium on Algebraic Geometry (Kinosaki, 1979), 9--23.

\bibitem[Mu81]{mu81}
S. Mukai, Duality between $D(X)$ and $D(\hat X)$ 
with its application to Picard sheaves, 
{\it Nagoya Math. J.}, {\bf 81} (1981), 153--175.

\bibitem[Ra78]{ra}
M. Raynaud, 
Contre-exemple au ``vanishing theorem'' en caract\'eristique $p>0$, 
C. P. Ramanujam --- A tribute, 
{\it Studies in Math.} {\bf 8} (1978), 273--278.

\bibitem[Sh94]{sha}
I. R. Shafarevich, 
{\it Basic algebraic geometry I}, second edition, 
Springer-Verlag, 1994.

\bibitem[SB91]{sb}
N. I. Shepherd-Barron, 
Geography for surfaces of general type in positive characteristic, 
{\it Inv. Math.}, {\bf 106} (1991), 263--274.

\bibitem[Ta72a]{ta72a}
H. Tango, 
On the behavior of extensions of vector bundles under the Frobenius map, 
{\it Nagoya Math. J.}, {\bf 48} (1972), 73--89.

\bibitem[Ta72b]{ta72b}
H. Tango, 
On the behavior of cohomology classes of vector bundles under 
the Frobenius map (in Japanese), Kyoto Univ., RIMS, Kokyuroku, 
{\bf 144} (1972), 93--102.

\bibitem[Xie05]{xie}
Q. Xie, On pseudo-effectivity of the second Chern classes for terminal 
threefolds, {\it Asian J. Math.}, {\bf 9} (2005), 121--132.

\bibitem[Xie]{xiep}
Q. Xie, Some remarks on the effective non-vanishing conjecture, 
preprint.

\end{thebibliography}
\end{document}